\documentclass[11pt,leqno]{article}
\usepackage{epic}
\usepackage{eepic}
\usepackage{amsfonts}
\newcommand{\achsbeschr}[2]{\makebox(0,0){$\frac{#1}{#2}$}}
\newcommand{\graphlinie}[2]{
        \put(#1,#2){\circle*{.1}}
        \put(#1,#2){\line(1,0){.9}}
                        }
\newcommand{\cN}{{\cal N}}

\newcommand{\Var}{\mbox{ Var}}
\newtheorem{theorem}{Theorem}[section]
\newtheorem{proposition}{Proposition}[section]

\newtheorem{corollary}{Corollary}[section]

\newtheorem{lemma}{Lemma}[section]

\def\Dto{\stackrel{{\cal D}}{\rightarrow}}
\def\Pto{\stackrel{{P}}{\rightarrow}}
\def\llog{\log\log}

\def\Deq{\stackrel{{\cal D}}{=}}

\def\tZ{\widetilde Z}

\newcommand{\beq}{\begin{equation}}
\newcommand{\eeq}{\end{equation}}
\newcommand{\beqa}{\begin{eqnarray}}
\newcommand{\eeqa}{\end{eqnarray}}
\newcommand{\beqaa}{\begin{eqnarray*}}
\newcommand{\eeqaa}{\end{eqnarray*}}
\newcommand{\bei}{\begin{itemize}}
\newcommand{\eei}{\end{itemize}}
\newcommand{\bee}{\begin{enumerate}}
\newcommand{\eee}{\end{enumerate}}

\begin{document}

\vspace*{1.5cm}

\noindent
{\Huge On Vervaat and Vervaat-error type processes\\
        for partial sums and renewals}

\vspace{1cm}
\noindent
{\large Endre Cs\'aki}\footnote{Research supported by the Hungarian
National Foundation for Scientific Research, Grant No. T 037886
and T 043037.}

\smallskip

\noindent
{\sl Alfr\'ed R\'enyi Institute of Mathematics, Hungarian Academy of
Sciences, Budapest, P.O.B. 127,
\newline H-1364, Hungary. }
{\texttt{~E-mail:~csaki@renyi.hu}}

\bigskip

\noindent
{\large Mikl\' os Cs\"org\H o}\footnote{Research supported by an NSERC
Canada  Grant at Carleton University, Ottawa.}

\smallskip

\noindent
{\sl School of Mathematics and Statistics, Carleton University,
Ottawa, Ontario, Canada K1S 5B6.}
\newline
\texttt{E-Mail:~mcsorgo@math.carleton.ca}

\bigskip

\noindent
{\large Zdzis\l aw Rychlik}\footnote{Research supported
by the Deutsche Forschungsgemeinschaft
through the German--Polish project 436 POL 113/98/0-1 "Probability measures".}

\smallskip\noindent
{\sl Instytut Matematyki, Uniwersytet Marii Curie-Sk\l odowskiej,
pl. Marii Curie-Sk\l odowskiej 1, \\
PL-20~031 Lublin, Poland. }
\texttt{E-Mail:~rychlik@golem.umcs.lublin.pl}

\bigskip

\noindent
{\large Josef Steinebach}\footnote{Research supported by a Polish-German
Exchange Grant No. BWM-III/DAAD/801/JK99, and by a special grant from Marie
Curie-Sk\l odowska University, Lublin.}

\smallskip\noindent
{\sl Universit\"at zu K\"oln, Mathematisches Institut, Weyertal 86-90,
D-50~931 K\"oln, Germany.}
\newline
\texttt{~E-Mail:~jost@math.uni-koeln.de}

\bigskip

\paragraph{Abstract.}
We study the asymptotic behaviour of stochastic processes that are
generated by sums of partial sums of i.i.d. random variables and
their renewals. We conclude that these processes cannot converge weakly to
any nondegenerate random element of the space $D[0,1]$. On the other hand
we show that their properly normalized integrals as Vervaat-type
stochastic processes converge weakly to a squared Wiener process.
Moreover, we also deal with the asymptotic behaviour of the deviations of
these processes, the so-called Vervaat-error type processes.




\section{Introduction}\label{sect1}

Let $X_1,X_2,\ldots$ be an i.i.d.\ sequence with d.f. $F(x) = P(X_1\leq x)$,
$-\infty <x<\infty$, and
\bei
\item[(i)]
$EX_1=\mu>0$;

\item[(ii)]
$0<\Var(X_1)=\sigma^2~(<\infty)$;

\item[(iii)]
$EX_1^4<\infty$.
\eei

\noindent
Set $S_n:=X_1+\cdots +X_n$ ($n=1,2,\ldots$), $S_0:=0$.  For $t\geq 0$
define pointwise the corresponding renewal (or first passage time) process
as $N(t):=\min\{n\geq 1~:~S_n>t\}$.
We further introduce the ``average'' and ``standardized'' processes
\beq\label{eq1.1}
S_n(t) := S_{[nt]}/(n\mu) \qquad (t\geq 0);
\eeq
\beq\label{eq1.2}
s_n(t) := n^{1/2} \sigma^{-1} \mu \Big( S_n(t)-t \Big)
=n^{-1/2}\sigma^{-1} \Big( S_{[nt]}-n\mu t \Big)
\qquad (t\geq 0);
\eeq
\beq\label{eq1.3}
N_n(t):=N(n \mu t)/n \qquad (t\geq 0);
\eeq
\beq\label{eq1.4}
r_n(t):= n^{1/2} \sigma^{-1} \mu \Big( N_n(t)-t \Big)
= n^{-1/2}\sigma^{-1}(\mu N(n\mu t)-n\mu t)
\qquad (t\geq 0),
\eeq
where $[x]$ denotes the integer part of $x$.

Throughout the paper $\Vert\cdot\Vert$ denotes the sup-norm on $D[0,1]$,
i.e., for $f\in D[0,1]$, $\Vert f\Vert:=\sup_{0\le t\le 1}|f(t)|$.

We first collect some well-known facts concerning the above processes.

\bigskip\noindent
{\bf Theorem A} {\it On a large enough probability space, for $X_1,
X_2,\dots$ satisfying {\rm (i), (ii)} and {\rm (iii)} above, one can
construct a standard Wiener process $\{W(t):~0\le t<\infty\}$, such that,
the processes $s_n$ and $r_n$ can be simultaneously approximated by $W$ as
follows:
\beq\label{eq1.5}
\lim_{n\to\infty} n^{1/4}\Vert s_n-W_n\Vert = 0 \quad \mbox{\rm a.s.}
\eeq
 and
\beq\label{eq1.6}
\limsup_{n\to\infty} n^{1/4}(\log n)^{-1/2}(\log\log n)^{-1/4}
\Vert r_n+W_n\Vert
= 2^{1/4} \sigma^{1/2} \mu^{-1/2} \quad \mbox{\rm a.s.},
\eeq
where $W_n(t):=n^{-1/2}W(nt)$.}

\bigskip
The result (\ref{eq1.5}) is due to Koml\'os {\it et al}. (1976), and
that of (\ref{eq1.6}) is due to Horv\'ath (1984). For more details and
further developments we refer to Chapter 2 of Cs\"org\H o and Horv\'ath
(1993).

Consider now the stochastic process
\beq\label{eq.1.6a}
R_n^*(t) := s_n(t) + r_n(t),\qquad 0\leq t\leq 1,
\eeq
which can be viewed as the remainder term in the Bahadur-Kiefer type
representation $r_n=-s_n+R_n^*$ of the renewal process $r_n$ in terms of the
partial sums process $s_n$. In particular, Bahadur (1966), Kiefer (1967,
1970) introduced and studied the stochastic process
$$
R_n(t):=\alpha_n(t)+\beta_n(t),\qquad 0\le t\le 1,
$$
that is known in the literature as the Bahadur--Kiefer process, where
\begin{eqnarray*}
    \alpha_n(t)
 &:=& n^{1/2}(F_n(t)-t),\qquad 0\le t\le 1,
    \\
    \beta_n(t)
 &:=& n^{1/2}(F_n^{-1}(t)-t),\qquad 0\le t\le 1,
\end{eqnarray*}
the uniform empirical and quantile processes, with $F_n^{-1}$ being the
left-continuous inverse of the right-continuously defined empirical
distribution function $F_n$ of the independent uniform (0,1) random
variables $U_1,\dots,U_n,\, n\ge 1$.

In this regard, we summarize the most relevant
results of Kiefer (1967, 1970) in the following theorem.

\medskip
\noindent
{\bf Theorem B} {\it For every fixed $t\in (0,1)$, we have
\begin{eqnarray}
    n^{1/4}R_n(t)
 &\Dto& (t(1-t))^{1/4} \cN (|{\widetilde \cN}|)^{1/2},\quad
    n\to\infty,
    \label{1.1} \\
    \limsup_{n\to\infty}{n^{1/4}| R_n(t)|\over (\log\log
n)^{3/4}}
 &=& (t(1-t))^{1/4}{2^{5/4}\over 3^{3/4}}\qquad {\rm a.s.},
    \label{1.3}
\end{eqnarray}
where $\cN$ and ${\widetilde \cN}$ are independent standard normal
variables and $\Dto$ denotes convergence in distribution. Also,}
\begin{equation}
   \lim_{n\to\infty} n^{1/4}(\log n)^{-1/2}{\| R_n\|\over
    (\|\alpha_n\|)^{1/2}}=1\qquad {\rm a.s.}
    \label{1.4}
\end{equation}

\medskip
As to (\ref{1.4}), Kiefer (1970) announced it but proved only convergence
in probability (cf. Theorem 1A, and the two sentences right after, in
Kiefer, 1970). The upper bound for the almost sure convergence in
(\ref{1.4}) was proved by Shorack (1982), and the lower bound by Deheuvels
and Mason (1990).

Concerning similar known results for the process $R_n^*$, we summarize
them in the next theorem.

\medskip\noindent
{\bf Theorem C} {\it Under the assumptions {\rm (i), (ii), (iii)} we
have for every fixed $t\in (0,1]$
\beq\label{eq1.10}
n^{1/4}R_n^*(t) \Dto t^{1/4} \sigma^{1/2} \mu^{-1/2}
\cN (|{\widetilde \cN}|)^{1/2},\quad
    n\to\infty.
\eeq
where, as in {\rm Theorem B}, $\cN$ and $\widetilde\cN$ are independent
standard normal random variables, and

\beq\label{eq1.7}
\lim_{n\to\infty} n^{1/4}(\log n)^{-1/2}\Vert R_n^*\Vert/\Vert r_n\Vert^{1/2}
= \sigma^{1/2} \mu^{-1/2} \quad \mbox{\rm a.s.}
\eeq
}

For the result in (\ref{eq1.10}) we refer to Cs\"org\H o and Horv\'ath
(1993, Theorem 2.1.5) and for that of (\ref{eq1.7}) to Deheuvels and Mason
(1990, Theorem 1B).

As a consequence of (\ref{eq1.7}), on account of

\beq\label{eq1.8}
\Vert r_n\Vert \Dto \Vert W\Vert \qquad \mbox{(cf.,\ e.g.,\ Vervaat
(1972))},
\eeq
as $n\to\infty$, we have (cf. Deheuvels and Mason (1990))

\beq\label{eq1.9}
n^{1/4}(\log n)^{-1/2}\Vert R_n^*\Vert \Dto \sigma^{1/2} \mu^{-1/2}
\Vert W\Vert.
\eeq

Based on (\ref{eq1.10}) and (\ref{eq1.9}), the next conclusion is
immediate.

\medskip\noindent
{\bf Corollary 1.1} {\it Given the assumptions {\rm (i), (ii), (iii)},
the statement
\beq\label{eq3.20}
a_nR_n^* \Dto Y, \quad  n\to\infty,
\eeq
cannot hold true in the space $D[0,1]$ (endowed with the Skorokhod
topology) for any non-degenerate random element $Y$ of $D[0,1]$ with any
normalizing sequence $\{a_n\}$.}

\bigskip

The respective results of (\ref{eq1.10}) and (\ref{eq1.7}) are based
on a strong invariance principle for the Bahadur-Kiefer type
process $\{R_n^*(t),\, 0\le t\le 1;\, n=1,2,\dots\}$ that was explicitly
stated in Cs\"org\H o and Horv\'ath (1993, p. 43) as follows.

\bigskip\noindent
{\bf Theorem D} {\it On assuming the conditions {\rm (i), (ii), (iii)},
on the probability space of {\rm Theorem A}, as $n\to\infty$, we have
\beq\label{invar}
R_n^*(t)=n^{-1/2}\left(W(nt)-W\left(nt-{\sigma\over
\mu}W(nt)\right)\right) +o(n^{-1/4})\qquad {\rm a.s.},
\eeq
uniformly in $t\in [0,1]$.}

\bigskip
Corollary 1.1 and Theorem 3 of Vervaat (1972) serve as motivation for
studying the following Vervaat-type process for partial sums and renewals:
\beqa\label{eq1.12}
V_n(t) &:=& \int^t_0 \Big\{ \Big(S_n(s)-s\Big) +
\Big(N_n(s)-s\Big)\Big\}ds\\
&=& \int^t_0 M_n(s)ds, \qquad 0\leq t \leq 1,\nonumber
\eeqa
where $M_n(s):=n^{-1/2}\sigma\mu^{-1}R_n^*(s),\, \, s\ge 0$.

\medskip
So far we have been dealing with renewal processes of partial sums under
the conditions (i), (ii) and (iii), that is to say we had general
renewal processes in mind. We will now see that the just introduced
Vervaat-type process $V_n$ will be well-behaving asymptotically as
$n\to\infty$. However the asymptotic behaviour of $V_n$ will be different
when it is based on general renewal processes, as compared
to it being based on ordinary renewal processes, i.e., when in addition
to (i), (ii), (iii), we also assume that $X_1$ is positive. Consequently,
in our Section \ref{sect2} we will deal with Vervaat-type processes for
ordinary renewals, while Section \ref{sect3} will be devoted to studying
such processes for general renewals. Moreover, in Section 4 we will also
be studying the asymptotic behaviour of Vervaat error-type processes that,
based on the respective results of Sections 2 and 3, will be defined there
(cf. (\ref{eq4.1})).

\setcounter{equation}{0}
\section{The Vervaat-type process for ordinary renewals}\label{sect2}

In addition to the conditions (i), (ii), (iii) of Section 1, in this
Section we assume also the condition

\medskip\noindent
(iv) $\quad P(X_1 > 0) = 1$.

\medskip\noindent

\begin{theorem}
Assume the conditions {\rm (i), (ii), (iii) and (iv)}. Then, as
$n\to\infty$
$$
{n\over \llog n} \Big\Vert V_n - {1\over 2} \bar S_n^2\Big\Vert \to 0 \quad
\mbox{\rm a.s.}
\leqno{\rm (a)}
$$
and
$$
n\Big\Vert V_n - {1\over 2} \bar S_n^2 \Big\Vert \Pto 0,
\leqno{\rm (b)}
$$
where $\bar S_n(t):=S_n(t)-t$.
\end{theorem}
\noindent
{\bf Proof}. Given condition (iv) it can be checked (cf. Figure~1 below)
that $V_n$ has the following representation:
\beq\label{eq2.1}
V_n(t) = A_n(t) - \bar S_n(t) \bar N_n(t) - {1\over 2} \bar N_n^2 (t),
\eeq
where
\beq\label{eq2.2}
A_n(t) := \int^t_{N_n(t)} \Big(\bar S_n(s) - \bar S_n(t)\Big)ds,
\eeq
\beq\label{eq2.3}
\bar S_n(t) := S_n(t)-t = (n\mu)^{-1}\Big(S_{[nt]}-n\mu t\Big),
\eeq
\beq\label{eq2.4}
\bar N_n(t) := N_n(t)-t = n^{-1}\Big(N(n\mu t) -n t\Big),
\eeq
for all $0\leq t\leq 1$.

In order to check the above representation, we note that
\[
N(n\mu t) = k, \qquad \hbox{for} \quad {S_{k-1}\over n\mu} \leq t <
{S_k\over n\mu}.
\]
Also, by calculations, we arrive at
\beqaa
\lefteqn{
\int^t_0 \Big(S_n(s)+N_n(s)\Big)ds = t^2 + \int^{N_n(t)}_t (t-S_n(s))ds}\\
&=& \int^t_{N_n(t)} \Big\{(S_n(s)-s)-(S_n(t)-t)\Big\}ds
- \int^t_{N_n(t)} (t-s)ds
-(S_n(t)-t)(N_n(t)-t) + t^2.
\eeqaa
Consequently, we conclude (\ref{eq2.1}) (cf. also Figure 1)
as follows:
\[
V_n(t)=\int^t_0 \Big( (S_n(s)-s) + (N_n(s)-s)\Big)ds
= A_n(t) - {1\over 2}\bar N_n^2(t) - \bar S_n(t) \bar N_n(t) + t^2-t^2.
\]



\begin{center}
\begin{picture}(14,11)(0,0)
\put(-.5,0){\vector(1,0){13.5}}
\multiput(1,0)(1,0){11}{\line(0,-1){.1}}

\multiputlist(.9,-.7)(1,0){
        \achsbeschr{1}{n}, \achsbeschr{2}{n},
         ,  ,  ,  ,
        \achsbeschr{[nt]}{n},
         ,  ,  ,
        \achsbeschr{N(n\mu t)}{n}
        }
\put(11.6,-.8){\makebox{$= N_n (t)$}}
\put(7.3,-.5){\makebox{$t$}}
\put(7.4,0){\line(0,-1){.1}}
\put(0,-.5){\vector(0,1){10.5}}

\put(-1,9){\achsbeschr{S_{N_n (t)}}{n\mu}}
\put(0,9){\line(-1,0){.1}}
\put(-.4,8.6){\makebox{$t$}}
\put(0,8.7){\line(-1,0){.1}}
\put(-1,8){\achsbeschr{S_{N_n (t)-1}}{n\mu}}
\put(0,8){\line(-1,0){.1}}
\put(-1,6){\achsbeschr{S_{[nt]}}{n\mu}}
\put(0,6){\line(-1,0){.1}}
\put(-1,3){\achsbeschr{S_{2}}{n\mu}}
\put(0,3){\line(-1,0){.1}}
\put(-1,2){\achsbeschr{S_{1}}{n\mu}}
\put(0,2){\line(-1,0){.1}}
\graphlinie{1}{2}
\graphlinie{2}{3}

\graphlinie{7}{6}
\graphlinie{8}{6.3}
\graphlinie{9}{7.1}
\graphlinie{10}{8}
\graphlinie{11}{9}

\thicklines
\dottedline{.2}(7.4,0)(7.4,9)
\dottedline{.2}(0,8.7)(12,8.7)
\dottedline{.2}(11,0)(11,9)
\dottedline{.2}(1,0)(1,2)
\dottedline{.2}(2,0)(2,3)
\dottedline{.2}(8,6)(8,6.3)
\dottedline{.2}(9,6.3)(9,7.1)
\dottedline{.2}(10,7.1)(10,8)
\dottedline{.2}(7.4,6.3)(7.7,6)
\dottedline{.2}(7.4,6.6)(8,6)
\dottedline{.2}(7.4,6.9)(8,6.3)
\dottedline{.2}(7.4,7.2)(8.3,6.3)
\dottedline{.2}(7.4,7.5)(8.6,6.3)
\dottedline{.2}(7.4,7.8)(8.9,6.3)
\dottedline{.2}(7.4,8.1)(9,6.5)
\dottedline{.2}(7.4,8.4)(9,6.8)
\dottedline{.2}(7.4,8.7)(9,7.1)

\dottedline{.2}(7.7,8.7)(9.3,7.1)
\dottedline{.2}(8,8.7)(9.6,7.1)
\dottedline{.2}(8.3,8.7)(9.9,7.1)
\dottedline{.2}(8.6,8.7)(10,7.3)
\dottedline{.2}(8.9,8.7)(10,7.6)
\dottedline{.2}(9.2,8.7)(10,7.9)
\dottedline{.2}(9.5,8.7)(10.2,8)
\dottedline{.2}(9.8,8.7)(10.5,8)
\dottedline{.2}(10.1,8.7)(10.8,8)
\dottedline{.2}(10.4,8.7)(11,8.1)
\dottedline{.2}(10.7,8.7)(11,8.4)
\put(8,8){\line(1,1){1.2}}
\put(9.3,9.3){$V_n (t)$}

\thinlines

\end{picture}

\end{center}


\bigskip
\bigskip
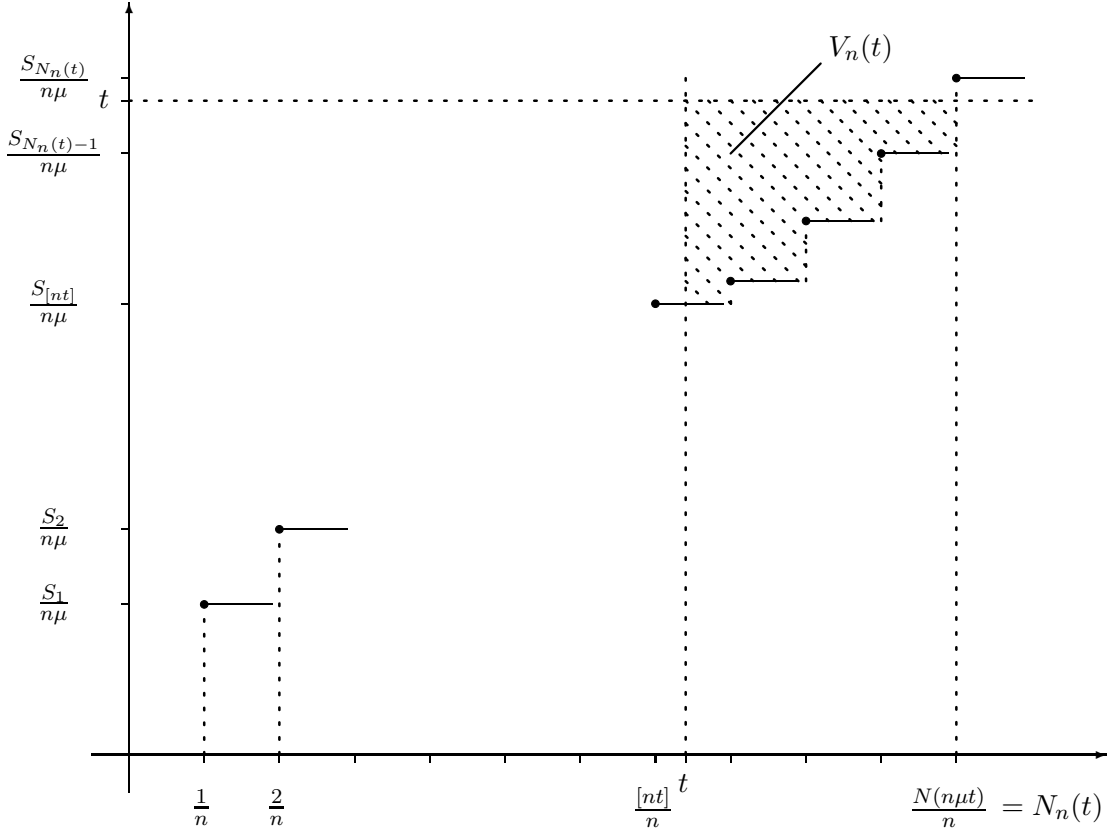
\centerline{Figure 1}
\vskip 1cm

As an immediate consequence in our proof, we get
\beqa\label{eq2.5}
V_n(t) - {1\over 2} \bar S_n^2(t) &=& A_n(t) - {1\over 2}\Big(\bar S_n(t) +
\bar N_n(t)\Big)^2\\
&=& A_n(t) - {1\over 2}M_n^2(t).\nonumber
\eeqa
Now, in view of (\ref{eq1.6}), (\ref{eq1.7}) and the law of the iterated
logarithm for a Wiener process,
\beq\label{eq2.6}
\Vert M^2_n\Vert = O\Big(n^{-3/2} (\log n)(\log \log n)^{1/2}\Big) \quad
\mbox{a.s.},
\eeq
and, as a consequence of (\ref{eq1.9}),
\beq\label{eq2.7}
\Vert M_n^2\Vert = O_P\Big(n^{-3/2} \log n\Big).
\eeq

Next, we show that $A_n$ is the dominating process on the right-hand
side of (\ref{eq2.5}). In (\ref{eq2.2}) put $s=t-u(t-N_n(t)) = t+u\bar
N_n(t)$. Then we get
\beq\label{eq2.8}
A_n(t) = - \bar N_n(t) \int^1_0\Big( \bar S_n(t+u \bar N_n(t))-\bar S_n(t)
\Big)du.
\eeq
On account of (\ref{eq1.5}),
\beq\label{eq2.9}
{n\mu\over\sigma} \bar S_n(t) = W(nt) + o\big(n^{1/4}\big) \quad \mbox{a.s.},
\eeq
uniformly in $t\in [0,1]$. A combination of (\ref{eq2.8}) and (\ref{eq2.9})
yields
\beq\label{eq2.10}
{n\mu\over\sigma} A_n(t) = - \bar N_n(t) \int^1_0
\Big(W(nt + un \bar N_n(t)) - W(nt)\Big)du
+ o\Big(n^{-1/4}(\log\log n)^{1/2}\Big) \quad \mbox{a.s.},
\eeq
uniformly in $t\in [0,1]$, since
\[
n \bar N_n(t) = O\Big( (n\log\log n)^{1/2}\Big) \quad \mbox{\rm a.s.},
\]
uniformly in $t\in [0,1]$.

This rate in (\ref{eq2.10}) can also be replaced by
$o_P(n^{-1/4})$,
on account of
\[
n\bar N_n(t) = O_P(n^{1/2}),
\]
uniformly in $t\in [0,1]$.

Moreover, by Theorem 1.2.1 of Cs\"org\H o and R\'ev\'esz (1981),
\beq\label{eq2.11}
W(nt + un\bar N_n(t)) - W(nt)
= O\Big(n^{1/4}(\log n)^{1/2}(\log\log n)^{1/4}\Big) \quad \mbox{\rm a.s.}
\eeq
as well as
\beq\label{eq2.12}
W(nt + un\bar N_n(t)) - W(nt)
= O_P\Big(n^{1/4}(\log n)^{1/2}\Big),
\eeq
uniformly in $u,t \in [0,1]$.

In view of (\ref{eq1.6}), we also have
\beq\label{eq2.13}
{\mu\over \sigma} n\bar N_n(t) =
-W(nt) + O\Big(n^{1/4}(\log n)^{1/2}(\log\log n)^{1/4}\Big) \quad
\mbox{\rm a.s.}
\eeq

As a consequence of (\ref{eq2.10})--(\ref{eq2.13}), we arrive at
\beq\label{eq2.14}
{n\mu\over \sigma} A_n(t) =
{\sigma W(nt)\over n\mu} \int^1_0
\Big(W(nt + un \bar N_n(t)) - W(nt)\Big)du
+ o\Big(n^{-1/4}(\log\log n)^{1/2}\Big) \quad \mbox{\rm a.s.},
\eeq
uniformly in $t\in [0,1]$, where this rate can also be replaced by
$o_P(n^{-1/4})$.

On making use of (\ref{eq2.13}) in combination with Theorem 1.2.1 of
Cs\"org\H o and R\'ev\'esz (1981), we conclude the following
strong and weak invariance principles for $A_n(t)$ of (\ref{eq2.2}):
\beqa\label{eq2.15}
~~~~A_n(t)&=& {\sigma W(nt)\over n\mu} \int^1_0 {\sigma\over n\mu}
\Big(W\Big(nt-u{\sigma\over\mu} W(nt)\Big)- W(nt)\Big)du
+ o\Big(n^{-5/4}(\log \log n)^{1/2}\Big) ~~ \mbox{a.s.}\\
&=& {\sigma\over n\mu} \int^{{\sigma\over n\mu}W(nt)}_0
(W(nt-nx)-W(nt))dx + o\Big(n^{-5/4}(\log\log n)^{1/2}\Big)\quad\mbox{a.s.}
\nonumber\\
&=& \int_0^{Y_n(t)} (Y_n(t-x)-Y_n(t))dx
+ o\Big(n^{-5/4}(\log\log n)^{1/2}\Big)\quad \mbox{a.s.}, \nonumber
\eeqa
uniformly in $t\in [0,1]$, where $Y_n(t) := {\sigma\over n\mu} W(nt)$, and
this a.s.\ rate can also be replaced by $o_P\Big(n^{-5/4}\Big)$.

Arguing similarly as in (\ref{eq2.11}) and (\ref{eq2.12}), we also have
\beqa\label{eq2.16}
A_n(t) &=&
O\Big(n^{-1/2}(\llog n)^{1/2}~ n^{-3/4}(\log n)^{1/2}(\llog n)^{1/4}\Big)
\quad \mbox{a.s.}\\
&=&
O\Big(n^{-5/4}(\log n)^{1/2}(\llog n)^{3/4}\Big) \quad \mbox{a.s.},\nonumber
\eeqa
and
\beq \label{eq2.17}
A_n(t) = O_P\Big( n^{-5/4}(\log n)^{1/2}\Big),
\eeq
uniformly in $t\in [0,1]$.

Thus, in view of (\ref{eq2.5}), (\ref{eq2.6}) and (\ref{eq2.16}) we
conclude (a) of Theorem 2.1, while (\ref{eq2.5}), (\ref{eq2.7}) and
(\ref{eq2.17}) result in (b) of Theorem 2.1.$\hfill\square$

As an immediate consequence of Theorem 2.1 and Theorem A, we get
\medskip\noindent
\begin{corollary}
Assume the conditions {\rm (i), (ii), (iii) and (iv)}. Then, as
$n\to\infty$, on the probability space of {\rm Theorem A} with $W_n$ as in
{\rm (\ref{eq1.5})}, we have
$$
{1\over \llog n} \Big\Vert nV_n - {1\over 2}
\left({\sigma\over\mu}W_n\right)^2\Big\Vert
\to 0 \quad
\mbox{\rm a.s.}
\leqno{\rm (a)}
$$
and
$$
\Big\Vert nV_n - {1\over 2} \left({\sigma\over\mu}W_n\right)^2
\Big\Vert \Pto 0.
\leqno{\rm (b)}
$$
\end{corollary}

As a consequence of (b) of Corollary 2.1, one concludes also
\beq\label{eq2.18}
nV_n \Dto Z, \quad n\to\infty,
\eeq
in the space $C[0,1]$ (endowed with the uniform topology), where $Z(t):=
{1\over 2}\Big({\sigma\over\mu} W(t)\Big)^2, \ t\in [0,1]$.

The strong invariance principle (a) of Corollary 2.1 in turn implies
Strassen (1964)-type laws of the iterated logarithm for $V_n$ as follows.

\begin{corollary}
Assume the conditions {\rm (i), (ii), (iii) and (iv)}. Then, as
$n\to\infty$, the set
\beq\label{flil}
\left\{{\mu^2 n V_n\over\sigma^2 \log\log n}, \, n\ge 3\right\}
\eeq
is relatively compact in $C[0,1]$, equipped with the sup-norm $\Vert
\cdot\Vert$. Furthermore, the set of all limit points of the functions
in the set {\rm (\ref{flil})} almost surely coincides with the set
$\{f^2:\, f\in{\cal S}\}$, where ${\cal S}$ is the Strassen class of all
absolutely continuous functions $f$ on $[0,1]$ such that $f(0)=0$ and
$\Vert f'\Vert_2\le 1$, where $\Vert\cdot\Vert_p$ denotes the $L_p$-norm,
$p\ge 1$.
\end{corollary}

As examples of consequences of Corollary 2.2 \`a la Strassen (1964), we
mention the following results.

\begin{corollary}
Under the conditions {\rm (i), (ii), (iii) and (iv)} we have
$$
\limsup_{n\to\infty}{n\Vert V_n\Vert \over \log\log n}=
{\sigma^2\over\mu^2}\qquad\mbox{\rm a.s.}
$$
and
$$
\limsup_{n\to\infty}{n\Vert V_n\Vert_1 \over \log\log n}=
{4\sigma^2\over\mu^2\pi^2}\qquad\mbox{\rm a.s.}
$$
\end{corollary}

\setcounter{equation}{0}
\section{The Vervaat process for general renewals}\label{sect3}

In this section we study the case of general renewals for which assumption
(iv) of Section~2 may not necessarily be true. In order to do this, we
need to introduce some further notations and auxiliary results:

\medskip\noindent
Let $\nu_i$ denote the $i$-th (strong) ascending ladder index of the partial
sum sequence $\{S_n\}_{n=0,1,\ldots}$, i.e., $\nu_0 = 0$, and,
recursively,
\beq\label{eq3.1}
\nu_i := \min\Big\{n>\nu_{i-1} \,:\, S_n-S_{\nu_{i-1}} > 0 \Big\}.
\eeq
We note that, under assumptions (i), (ii), (iii),
$\{\nu_i-\nu_{i-1}\}_{i=1,2,\ldots}$
is an i.i.d.\ sequence of random variables with
\beq\label{eq3.2}
E\nu_1^4 < \infty.
\eeq
Moreover, the sequence $\{ S_{\nu_i} - S_{\nu_{i-1}} \}_{i=1,2,\ldots}$ of
ladder heights is also an i.i.d.\ sequence of random variables with
\beq\label{eq3.3}
ES^4_{\nu_1} < \infty
\eeq
(cf., e.g., Gut (1988, Sections III. 2-3)).  By definition,
\[
N(t) = \nu_i, \quad \mbox{for } S_{\nu_{i-1}} \leq t<S_{\nu_i} \quad (i=1,2,
\ldots; \, t\geq 0),
\]
i.e.,
\beq\label{eq3.4}
N(n\mu s) = \nu_i, \quad \mbox{for }  {S_{\nu_{i-1}}\over n\mu} \leq s <
{S_{\nu_i}\over n\mu}
\quad (i=1,2,\ldots; \, s\geq 0).
\eeq

Let $N(n\mu t) = \nu_\ell$, i.e.,
\beq\label{eq3.5}
\ell = \min\Big\{ j : S_{\nu_j} > n\mu t\Big\},
\eeq
that is, $\ell=N_H(n\mu t)$, where $N_H$ denotes the (ordinary) renewal
process corresponding to the (ladder height) sequence $S_{\nu_1},
S_{\nu_2} - S_{\nu_1},\ldots$.  This will play a crucial role in the
calculations below.

We first note that we continue to use the same definition for the Vervaat
process $V_n$ for general renewals as in (\ref{eq1.12}). Namely, we study
the integral
$$
V_n(t)=\int_0^t M_n(s)\, ds = \int^t_0\Big(\bar S_n(s) + \bar
N_n(s)\Big)ds,
$$
where $\bar S_n$ and $\bar N_n$ are as in (\ref{eq2.3}) and (\ref{eq2.4}),
respectively, i.e.,
$$
\bar S_n(t)=S_n(t)-t=(n\mu)^{-1}(S_{[nt]}-n\mu t),
$$
$$
\bar N_n(t)=N_n(t)-t=n^{-1} (N(n\mu t)-n t),
$$
with $S_n(t)$ and $N_n(t)$ as in (\ref{eq1.1}) and (\ref{eq1.3}),
respectively.

Consider now the following figure.



\begin{center}
\begin{picture}(15,11)(0,0)
\put(-1,0){\vector(1,0){14}}

\multiput(1,0)(1,0){11}{\line(0,-1){.1}}
\put(1,0){\line(0,1){.1}}
\put(2,0){\line(0,1){.1}}

\multiputlist(.9,-.6)(1,0){
 { },
 { },
        \achsbeschr{3}{n},
         ,
        \achsbeschr{\nu_2}{n},
         ,
        \achsbeschr{[nt]}{n},
         ,
        \achsbeschr{\nu_{\ell -1}}{n},
         ,
        \achsbeschr{N(n\mu t)}{n}
        }
\put(1,0.6){\achsbeschr{1}{n}}
\put(2,0.6){\achsbeschr{\nu_1}{n}}

\put(11.5,-.7){$= N_n (t) = \frac{\nu_{\ell}}{n}$}
\put(7.3,-.4){\makebox{$t$}}
\put(7.4,0){\line(0,-1){.1}}
\put(0,-.4){\vector(0,1){10}}  

\put(-1,8){\achsbeschr{S_{\nu_{\ell}}}{n\mu}}
\put(0,8){\line(-1,0){.1}}
\put(-.4,6.9){\makebox{$t$}}
\put(0,7){\line(-1,0){.1}}
\put(-1,6){\achsbeschr{S_{\nu_{\ell -1}}}{n\mu}}
\put(0,6){\line(-1,0){.1}}
\put(-1,3){\achsbeschr{S_{\nu_2}}{n\mu}}
\put(0,3){\line(-1,0){.1}}
\put(-1,2){\achsbeschr{S_{\nu_1}}{n\mu}}
\put(0,2){\line(-1,0){.1}}
\graphlinie{1}{-1}
\graphlinie{2}{2}
\graphlinie{3}{.7}
\graphlinie{4}{1.5}
\graphlinie{5}{3}

\graphlinie{7}{4.8}
\graphlinie{8}{4.3}
\graphlinie{9}{6}
\graphlinie{10}{5}
\graphlinie{11}{8}

\thicklines
\dottedline{.2}(7.4,0)(7.4,7.5)
\dottedline{.2}(0,7)(11.5,7)
\dottedline{.2}(11,0)(11,8)
\dottedline{.3}(1,0)(1,-1)
\dottedline{.3}(2,0)(2,-1)

\dottedline{.2}(1,-.8)(1.2,-1)
\dottedline{.2}(1,-.5)(1.5,-1)
\dottedline{.2}(1,-.2)(1.8,-1)
\dottedline{.2}(1.1,0)(2,-.9)
\dottedline{.2}(1.4,0)(2,-.6)
\dottedline{.2}(1.7,0)(2,-.3)

\put(.7,-.5){\line(1,0){.7}}
\put(.2,-.7){$A_1$}
\dottedline{.2}(3,2)(5,2)
\dottedline{.2}(3,.7)(3,2)
\dottedline{.2}(4,.7)(4,1.5)
\dottedline{.2}(5,1.5)(5,3)

\dottedline{.2}(3,1)(3.3,.7)
\dottedline{.2}(3,1.3)(3.6,.7)
\dottedline{.2}(3,1.6)(3.9,.7)
\dottedline{.2}(3,1.9)(4,.9)
\dottedline{.2}(3.2,2)(4,1.2)
\dottedline{.2}(3.5,2)(4,1.5)
\dottedline{.2}(3.8,2)(4.3,1.5)
\dottedline{.2}(4.1,2)(4.6,1.5)
\dottedline{.2}(4.4,2)(4.9,1.5)
\dottedline{.2}(4.7,2)(5,1.8)

\put(3.6,1.7){\line(1,1){.5}}
\put(4.2,2.3){$A_2$}
\dottedline{.1}(8,4.8)(9,4.8)
\dottedline{.2}(8,4.3)(8,4.8)
\dottedline{.2}(9,4.3)(9,6)

\dottedline{.2}(8,4.6)(8.3,4.3)
\dottedline{.2}(8.1,4.8)(8.6,4.3)
\dottedline{.2}(8.4,4.8)(8.9,4.3)
\dottedline{.2}(8.7,4.8)(9,4.5)

\put(8.5,3.9){\line(0,1){.6}}
\put(8.2,3.4){$A_{\ell -1}$}

\dottedline{.1}(10,6)(11,6)
\dottedline{.2}(10,5)(10,6)

\dottedline{.2}(10,5.3)(10.3,5)
\dottedline{.2}(10,5.6)(10.6,5)
\dottedline{.2}(10,5.9)(10.9,5)
\dottedline{.2}(10.2,6)(11,5.2)
\dottedline{.2}(10.5,6)(11,5.5)
\dottedline{.2}(10.8,6)(11,5.8)

\put(10.5,4.6){\line(0,1){.6}}
\put(10.2,4.1){$A_{\ell}$}

\dottedline{.2}(7.4,6.8)(7.6,7)
\dottedline{.2}(7.4,6.4)(8,7)
\dottedline{.2}(7.4,6)(8.4,7)
\dottedline{.2}(7.4,5.6)(8.8,7)
\dottedline{.2}(7.4,5.2)(9.2,7)
\dottedline{.2}(7.4,4.8)(9.6,7)

\dottedline{.2}(7.8,4.8)(10,7)
\dottedline{.2}(8.2,4.8)(9,5.6)   \dottedline{.2}(9.4,6)(10.4,7)
\dottedline{.2}(8.6,4.8)(9,5.2)   \dottedline{.2}(9.8,6)(10.8,7)
                                  \dottedline{.2}(10.2,6)(11,6.8)
                                  \dottedline{.2}(10.6,6)(11,6.4)
                                  \dottedline{.2}(11,6)(11,6)

\thinlines

\end{picture}
\end{center}

\bigskip\bigskip
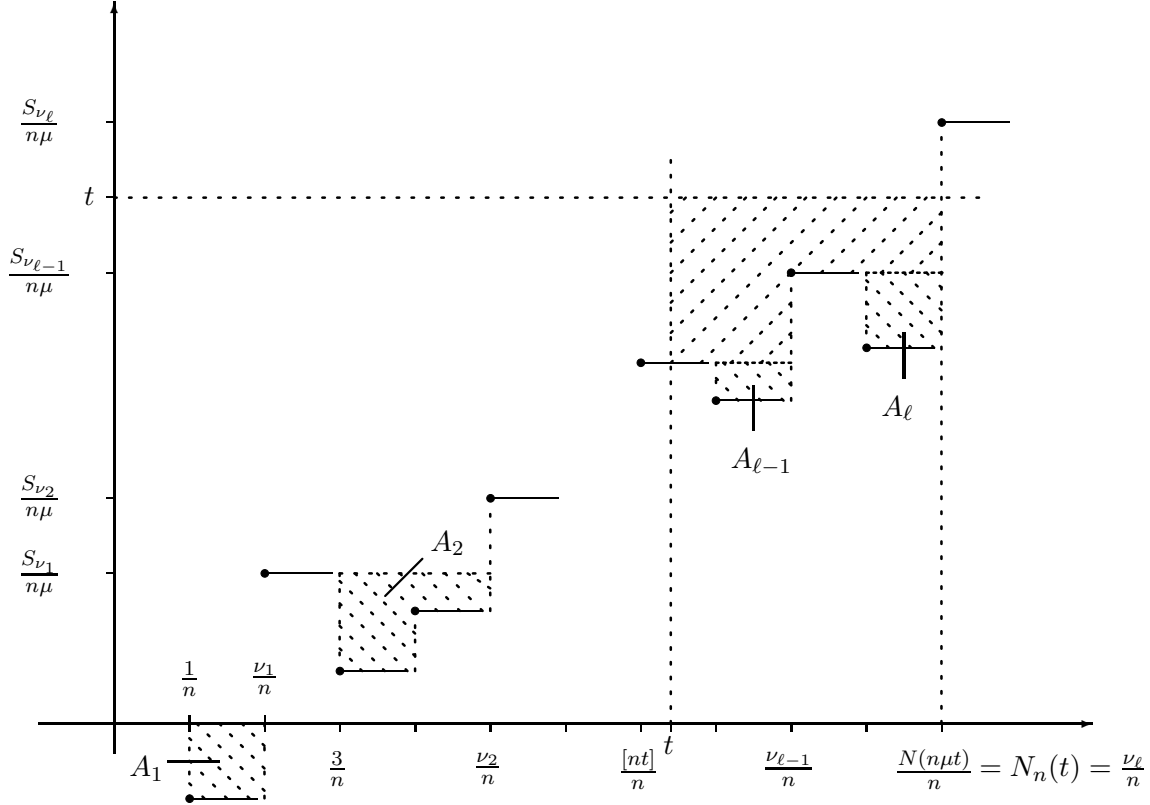
\centerline{Figure 2}

\vskip 1cm

Note that, in contrast to the case of the ordinary renewal process (cf.
Figure 1), the (random) areas $A_1,A_2,\ldots$ (cf. Figure~2) are
additional in the present case of the Vervaat process for general
renewals. Thus, in computing the above integral $V_n$, these areas are to
be taken into consideration. Our next lemma accomplishes this via
establishing the following new representation for $V_n$ in the general
case that will now replace that of (\ref{eq2.5}).

\begin{lemma}
Under the conditions {\rm (i), (ii), (iii)} we have
\beq\label{eq3.16}
V_n(t)= A_n(t) + B_n(t) - {1\over 2} M_n^2(t)+{1\over 2} \bar S_n^2(t),
\eeq
with $A_n$ as in {\rm (\ref{eq2.2})}, $M_n=\bar S_n+\bar N_n$ and $B_n$ is
defined by
\beq\label{eq3.7}
B_n(t) := {1\over n^2\mu} \sum^{N_H(n\mu t)}_{i=1}
\left( \sum^{\nu_i-\nu_{i-1}-1}_{j=1}(\nu_i-\nu_{i-1}-j)X_{\nu_{i-1}+j}
\right),
\eeq
with $N_H(n\mu t):=\ell$ of {\rm (\ref{eq3.5})}.
\end{lemma}

\noindent{\bf Proof}. We first show that
\beq\label{eq3.6}
\int^t_0 \bar N_n(s)ds + \int^{N_n(t)}_0 \bar S_n(s)ds = -{1\over 2}
\bar N_n^2(t) + B_n(t),
\eeq
We have
\beq\label{eq3.8}
\int^t_0 (N_n(s)-s)ds + \int^{N_n(t)}_0(S_n(s)-s)ds
= \int^t_0 N_n(s)ds + \int^{N_n(t)}_0 S_n(s)ds - {t^2\over 2}-{N_n^2(t)\over 2}.
\eeq
By (\ref{eq3.4}) and the definition of $N_n$ we arrive at
\begin{eqnarray}\label{3.9}
\int^t_0 N_n(s)ds
&=& {1\over n} \sum^\ell_{i=1}
\int^{{S_{\nu_i}\over n\mu}}_{{S_{\nu_{i-1}}\over n\mu}} \nu_i~ds-{1\over n}
\int_t^{{S_{\nu_\ell}\over n\mu}} \nu_\ell~ds\\
&=& {1\over n^2\mu} \sum^\ell_{i=1}
\nu_i \Big(S_{\nu_i}-S_{\nu_{i-1}}\Big) -
{1\over n} \nu_\ell \Big( {S_{\nu_\ell}\over n\mu} - t \Big).\nonumber
\end{eqnarray}
Also,
\begin{eqnarray}\label{3.10}
\int^{N_n(t)}_0 S_n(s)ds
&=&
\sum^{\nu_\ell}_{i=1} \int^{i\over n}_{i-1\over n} S_n(s)ds
= {1\over n^2\mu} \sum^{\nu_\ell}_{i=1} S_{i-1}\\
&=&{1\over n^2\mu} \sum^\ell_{i=1} \sum^{\nu_i-\nu_{i-1}}_{j=1}
\Big(S_{\nu_{i-1}+j-1} - S_{\nu_{i-1}} \Big)
+ {1\over n^2\mu} \sum^\ell_{i=1}(\nu_i - \nu_{i-1})S_{\nu_{i-1}}.\nonumber
\eeqa

Now, for any real sequences $\{a_i\}_{i=0,1,\ldots}$ and $\{b_i\}_
{i=0,1,\ldots}$, it can easily be checked that
\[
\sum^n_{i=1} a_i(b_i-b_{i-1}) + \sum^n_{i=1} b_{i-1}(a_i-a_{i-1}) = a_nb_n -
a_0b_0.
\]
Hence, since $\nu_0 = 0 = S_{\nu_0}$,
\beq\label{eq3.11}
\sum^\ell_{i=1} \nu_i\Big(S_{\nu_i} - S_{\nu_{i-1}}\Big)
= \nu_\ell S_{\nu_\ell} - \sum^\ell_{i=1}(\nu_i-\nu_{i-1})S_{\nu_{i-1}}.
\eeq
Thus, on combining (\ref{eq3.8}) - (\ref{eq3.11}), we obtain
\beqa\label{eq3.12}
\lefteqn{\int^t_0 \bar N_n(s)ds + \int^{N_n(t)}_0 \bar S_n(s) ds}\\
&=& {1\over n^2\mu} \sum^\ell_{i=1} \sum^{\nu_i-\nu_{i-1}-1}_{j=1}
\Big(S_{\nu_{i-1}+j} - S_{\nu_{i-1}}\Big)
- {t^2\over 2} - {N^2_n(t)\over 2} + t {\nu_\ell\over n}\nonumber\\
&=& - {1\over 2}(N_n(t) - t)^2 + {1\over n^2\mu}
\sum^\ell_{i=1}
\left(\sum^{\nu_i-\nu_{i-1}-1}_{j=1}
(\nu_i-\nu_{i-1}-j)X_{\nu_{i-1}+j} \right).\nonumber
\eeqa
This, in turn, proves (\ref{eq3.6}).

As a consequence of (\ref{eq3.6}) and the definition of $A_n(t)$ we
now conclude our claim in (\ref{eq3.16}) as follows.

\beqaa
V_n(t) &=& \int^t_0\Big(\bar S_n(s) + \bar N_n(s)\Big)ds\\
&=& \int^t_0 \bar N_n(s)ds + \int^{N_n(t)}_0 \bar S_n(s)ds
+ \int^t_{N_n(t)} \bar S_n(s)ds\\
&=& - {1\over 2} \bar N_n^2(t) + B_n(t) + A_n(t) - \bar S_n(t) \bar N_n(t)\\
&=& A_n(t) + B_n(t) - {1\over 2} M_n^2(t) + {1\over 2} \bar S_n^2(t).
\eeqaa
$\hfill\square$

For the sake of stating and proving an analogue of Theorem 2.1 in case of
the general renewal process, we need further notations.

Consider the summands

\beq\label{eq3.13}
D_i := \sum^{\nu_i-\nu_{i-1}}_{j=2} (\nu_i-\nu_{i-1}-j)X_{\nu_{i-1}+j-1}
\quad i=1,2,\ldots
\eeq
in the definition of $B_n$ (cf. \ref{eq3.7}). We note that the just
introduced random variables $D_i$ are i.i.d. with mean $\mu_D:=ED_1$
and finite second moment $ED_1^2<\infty$. The latter is implied by
(\ref{eq3.2}) and (\ref{eq3.3}), on applying the Cauchy-Schwarz inequality
when computing $ED_1^2$. Put also $\mu_H:=ES_{\nu_1}$, the
so-called expected ladder height (cf.,\ e.g.,\ Gut (1988, Chap.\ IV,
Theorem 2.7)).

The extension of our Theorem 2.1 to the general renewal
situation now reads as follows.

\begin{theorem}
Assume the conditions {\rm (i), (ii), (iii)}.  Then, as $n\to\infty$,
$$
{1\over\llog n} \Big|\Big| nV_n - {n\over 2} \bar S_n^2 - \Big({\mu_D\over
\mu_H}\Big)id \Big|\Big|
\to 0 \quad \mbox{\rm a.s.}
\leqno{\rm (a)}
$$
$$
\Big|\Big| nV_n - {n\over 2} \bar S_n^2 - \Big({\mu_D\over \mu_H}\Big)id
\Big|\Big|
\Pto 0~,
\leqno{\rm (b)}
$$
where $id(t):=t$, $t\in [0,1]$.
\end{theorem}

\medskip\noindent{\bf Proof}. For $B_n$ of Lemma 3.1 we conclude the
following asymptotic representation. As $n\to\infty$, by the law of the
iterated logarithm, we have almost surely

\beqa\label{eq3.17}
B_n(t) &=& {1\over n^2\mu} \sum^{N_H (n\mu t)}_{i=1} D_i\\
&=& {1\over n^2\mu}
\Bigg(\sum^{N_H(n\mu t)}_{i=1} D_i-n\mu t{\mu_D\over \mu_H} \Bigg)
+ {1\over n}\Big( {\mu_D\over \mu_H} \Big)t\nonumber\\
&=& {1\over n}\Big({\mu_D\over \mu_H}\Big)t
+ O\Big( {(\llog n)^{1/2}\over n^{3/2} } \Big),\nonumber
\eeqa
uniformly in $t\in [0,1]$. Consequently, via (\ref{eq2.6})/(\ref{eq2.7}),
(\ref{eq2.16})/(\ref{eq2.17}), Lemma 3.1 and (\ref{eq3.17}) we conclude
the respective two statements of Theorem 3.1.$\hfill\square$

\medskip \noindent
{\bf Remark 3.1} Via the definition of the ladder index $\nu_1$, we have
\[
X_1, X_1+X_2,\ldots,X_1 + \cdots +X_{\nu_1-1}\leq 0, ~X_1+\cdots + X_{\nu_1} >
0.
\]

Hence (cf. (\ref{eq3.13})) we have also
\[
D_1=\sum^{\nu_1-1}_{j=1} (\nu_1-j)X_j \leq 0.
\]

If, however,
$$
P(X_1 \geq 0) = 1,
\leqno{\rm (iv')}
$$
then
\beq\label{eq3.14}
P(D_1=0) = 1.
\eeq
On the other hand, if
$$
P(X_1<0)>0,
\leqno{\rm (v)}
$$
then we have
\beq\label{eq3.15}
P(D_1<0)>0,
\eeq
as well.

This shows that the results of Section~2 also remain true under assumption
(iv') replacing (iv), for then $B_n=0$ a.s.

As an immediate consequence of Theorem 3.1 and Theorem A, we conclude

\medskip\noindent{\bf Corollary 3.1}.
{\it Assume the conditions {\rm (i), (ii), (iii)}. Then, as $n\to\infty$,
on the probability space of {\rm Theorem A} with $W_n$ as in {\rm
(\ref{eq1.5})}, we have
$$
{1\over\llog n} \Big|\Big| nV_n - {1\over
2}\left({\sigma\over\mu}W_n\right)^2 - \Big({\mu_D\over
\mu_H}\Big)id \Big|\Big|
\to 0 \quad \mbox{\rm a.s.}
\leqno{\rm (a)}
$$
and
$$
\Big\Vert nV_n - {1\over 2}\Big({\sigma\over\mu} W_n\Big)^2 -
\Big({\mu_D \over \mu_H}\Big)id\Big\Vert \Pto 0,
\leqno{\rm (b)}
$$
}

As a consequence of (b) of Corollary 3.1, we conclude also
\beq\label{eq3.19}
nV_n \Dto \tZ, \quad n\to\infty,
\eeq
in the space $C[0,1]$ (endowed with the uniform topology),
where $\tZ(t):={1\over 2}\Big({\sigma\over\mu} W(t)\Big)^2 +
\Big({\mu_D\over \mu_H} \Big)t$, $t\in [0,1]$.

The latter conclusion in turn also allows us a Vervaat (1972)-type
proof of Corollary 1.1. Namely, if the statement (\ref{eq3.20}) were true,
then the sequence in there could be replaced by $\{n\}$, due to the
convergence of types theorem (cf. Lemma 1 on p. 253 of Feller (1971)). But
then one would also have the equality in distribution

\[
\int^t_0 Y(s)ds ~\Deq~ {1\over 2}\bigg({\sigma\over\mu}
W(t)\bigg)^2+\left({\mu_D\over\mu_H}\right) t,\qquad t\in [0,1],
\]
which is impossible, since the Wiener process $W$ is a.s.\ nowhere
differentiable. Consequently, the relation (\ref{eq3.20}) cannot hold
true.

In concluding this section, it may be of interest to note that, while it
is true that the asymptotic behaviour of the Bahadur-Kiefer-type
process $R_n^*$ of (\ref{eq1.7}) with ordinary and general renewals
is the same pointwise and in sup-norm, their integrated form $V_n$ behaves
differently as $n\to\infty$. Namely, in Corollary 3.1 and (\ref{eq3.19})
the squared Wiener process appears with a drift term which is not present
in Corollary 2.1 and (\ref{eq2.18}).

\bigskip
\setcounter{equation}{0}
\section{Vervaat-error type processes for partial sums and
renewals}\label{sect4}

\bigskip
Theorem 2.1 and Theorem 3.1 suggest the study of the asymptotic
behaviour of the following Vervaat-error type process:

\beq\label{eq4.1}
Q_n(t):=V_n(t)-{1\over 2}\bar S_n^2(t)-\left({\mu_D\over\mu_H}\right)t,
\qquad 0\le t\le 1.
\eeq

We note that the drift term $(\mu_D/\mu_H) t$ is identically zero in case
of an ordinary renewal process (cf. Theorem 2.1).

Define the stochastic process $Z_n$ by

\beq\label{eq4.2}
Z_n(t) := \int^{Y_n(t)}_0 \Big(Y_n(t-x)-Y_n(t)\Big)dx, \qquad 0\le t\le 1,
\eeq
where $Y_n(t) := {\sigma\over n\mu} W(nt)$.

We have the following representations for $Q_n$ in terms of $Z_n$.

\begin{proposition} Assume the conditions {\rm (i), (ii), (iii)}. Then, as
$n\to\infty$, on the probability space of {\rm Theorem A}, we have
\beq\label{eq4.3}
Q_n(t) = Z_n(t)+o\Big(n^{-5/4}(\llog n)^{1/2}\Big) \quad \mbox{\rm a.s.},
\eeq
and
\beq\label{eq4.4}
Q_n(t) = Z_n(t) + o_P\Big(n^{-5/4}\Big),
\eeq
uniformly in $t\in [0,1]$.
\end{proposition}

\medskip\noindent
{\bf Proof}. In view of (\ref{eq3.16}), via (\ref{eq2.15}),
(\ref{eq3.17}), and (\ref{eq2.6}), we conclude (\ref{eq4.3}). Similarly,
and again in view of (\ref{eq3.16}), via (\ref{eq2.15}) with the error term
$o_P(n^{-5/4})$, (\ref{eq3.17}), and (\ref{eq2.7}), we arrive at
(\ref{eq4.4}) as well.

\begin{theorem}
Assume {\rm (i), (ii), (iii)}.  Then, for $0<t\le 1$ fixed, as
$n\to\infty$
\beq\label{eq4.5}
n^{5/4} Q_n(t) \Dto  \Big({1\over 3}\Big)^{1/2}
\Big({\sigma\over\mu}\Big)^
{5/2} t^{3/4}\cN|{\widetilde \cN}|^{3/2},
\eeq
where $\cN$ and ${\widetilde \cN}$ are independent standard normal
random variables.
\end{theorem}

We note in passing that the result of (\ref{eq4.5}) parallels that of
(2.5) of Theorem 2.1 of Cs\'aki {\it et al}. (2002), which deals with a
similar problem concerning the Vervaat-error process of the empirical and
quantile processes. For further related results along these lines we refer
to Cs\"org\H o and Zitikis (1999).

\medskip
\noindent
{\bf Proof}. Given (\ref{eq4.4}) of Proposition 4.1, it is clear that the
pointwise asymptotic distribution of $Q_n(t)$ for any fixed $t\in (0,1]$
can be established via that of $Z_n(t)$. Clearly, from (\ref{eq4.2}) we
obtain
\beq\label{eq4.6}
Z_n(t) = {\sigma W(nt)\over n\mu} \int^1_0 {\sigma\over n\mu}
\Big(W(nt - s{\sigma\over\mu} W(nt))-W(nt)\Big)ds.
\eeq

Let
\beqaa
W_1(u) &:=& (nt)^{-1/2}(W(nt-unt)-(1-u)W(nt))+u\cN^*, \quad 0\leq u\leq
1,\\
W_2(u) &:=& (nt)^{-1/2}(W(nt+unt)-W(nt)), \quad 0\leq u\leq 1,
\eeqaa
where $\cN^*$ is a standard normal random variable, independent of
$W(\cdot)$. It is obvious that for any fixed $n$ and $t$, the
processes $\{W_1(u),\, 0\le u\le 1\}$ and $\{W_2(u),\, 0\le u\le 1\}$ are
independent standard Wiener processes, which are also independent of
$W(nt)$.

Consider $Z_n(t)$ of (\ref{eq4.6}), and first assume that $W(nt)>0$. Put
$unt=s(\sigma/\mu)W(nt)$ in (\ref{eq4.6}). Thus we obtain
$$
Z_n(t)={t\sigma\over n\mu}\int_0^{\sigma W(nt)\over\mu nt}
(W(nt-unt)-W(nt))\, du
$$
$$
={t\sigma\over\mu}\left({t\over n}\right)^{1/2}
\int_0^{\sigma W(nt)\over\mu nt} W_1(u)\, du
-{t\sigma\over\mu n}\int_0^{\sigma W(nt)\over\mu nt} uW(nt)\, du
-{t\sigma\over\mu}\left({t\over n}\right)^{1/2}
\int_0^{\sigma W(nt)\over\mu nt} u\cN^*\, du
$$
$$
={t\sigma\over\mu}\left({t\over n}\right)^{1/2}
\int_0^{\sigma W(nt)\over\mu nt} W_1(u)\, du +O_P(n^{-3/2}).
$$

Assuming now that $W(nt)<0$, on putting $unt=-s(\sigma/\mu)W(nt)$ in
(\ref{eq4.6}), we get
$$
Z_n(t)={t\sigma\over\mu}\left({t\over n}\right)^{1/2}
\int_0^{-\sigma W(nt)\over\mu nt} W_2(u)\, du.
$$

Hence, for each fixed $t\in (0,1]$,  $Z_n(t)$ has the following
representation
\beq\label{eq4.7}
Z_n(t)=Z_n^*(t)+O_P(n^{-3/2}),
\eeq
where
$$
Z_n^*(t):={t\sigma\over\mu}\left({t\over n}\right)^{1/2}
\int_0^{\left | \sigma W(nt)\over\mu nt\right|} W^*(u)\, du,
$$
and $W^*(u):= W_1(u)I\{W(nt)>0\}+W_2(u)I\{W(nt)<0\}$ is a standard
Wiener process, independent of $W(nt)$ for each fixed $n$ and $t$.

It is well-known that for any fixed $T>0$,
\[
\int^T_0 W^*(u)\, du ~\Deq~ N\Big(0,{T^3\over 3}\Big) ~\Deq~ {T^{3/2}\over
3^{1/2}} N(0,1).
\]
Consequently, by independence of $W(nt)$ and $\{W^*(u),\, 0\le u\le 1\}$
for each fixed $n$ and $t$,
\beq
Z_n^*(t) \Deq
\Big({1\over 3}\Big)^{1/2}\Big({\sigma\over\mu}\Big)^{5/2} {1\over n^2}
\bigg({|W(nt)|\over (nt)^{1/2}}\bigg)^{3/2} (nt)^{3/4}\cN.
\eeq
Hence for each fixed $n$ and $t$ we conclude
\beq
Z_n^*(t) \Deq
 n^{-5/4}\Big({1\over 3}\Big)^{1/2}\Big({\sigma\over\mu}\Big)^{5/2}
t^{3/4}|{\widetilde \cN}|^{3/2}\cN,
\eeq
where $\cN$ and ${\widetilde \cN}$ are independent standard normal random
variables.

Consequently, on account of (\ref{eq4.7}), as $n\to\infty$, we obtain
\beq\label{eq4.8}
n^{5/4} Z_n(t) \Dto  \Big({1\over 3}\Big)^{1/2}
\Big({\sigma\over\mu}\Big)^
{5/2} t^{3/4}\cN|{\widetilde \cN}|^{3/2},
\eeq
for each fixed $t\in(0,1]$.

This in view of (\ref{eq4.4}) of Proposition 4.1 also completes the proof
of (\ref{eq4.5}) of Theorem 4.1.$\hfill\square$



\bigskip


\end{document}